\documentclass[12pt]{amsart}
\usepackage{amsfonts,amsmath,amssymb,amsthm}
\usepackage{latexsym}
\usepackage{mathrsfs}
\usepackage{tikz}
\usepackage{color}

\newtheorem{lem}{Lemma}[section]
\newtheorem{prop}{Proposition}[section]
\newtheorem{thm}{Theorem}[section]
\newtheorem{rem}{Remark}
\newtheorem{corollary}{Corollary}[section]

\newcommand{\su}{\mathfrak{sl}}

\newcommand{\id}{\text{id}}
\newcommand{\cR}{\mathcal{R}}

\newcommand{\mC}{\mathfrak{C}}

\begin{document}

\title[Revisiting Askey--Wilson algebra]{Revisiting the Askey--Wilson algebra \\ with the universal $R$-matrix of $U_q(\su_2)$}
%\dedicatory{}

\author[N.Cramp\'e]{Nicolas Cramp\'e$^{\dagger}$}
\address{$^\dagger$ Institut Denis-Poisson CNRS/UMR 7013 - Universit\'e de Tours - Universit\'e d'Orl\'eans, 
Parc de Grammont, 37200 Tours, France.}
\email{crampe1977@gmail.com}

\author[J.Gaboriaud]{Julien Gaboriaud$^{*}$}
\address{$^*$ Centre de recherches math\'ematiques, Universit\'e de Montr\'eal,
P.O. Box 6128, Centre-ville Station,
Montr\'eal (Qu\'ebec), H3C 3J7, Canada.}

\author[L.Vinet]{Luc Vinet$^{*}$}

\author[M.Zaimi]{Meri Zaimi$^{*}$}

\begin{abstract}
A description of the embedding of a centrally extended  Askey--Wilson algebra, $AW(3)$, in $U_q(\su_2)^{\otimes 3}$ is given 
in terms of the universal $R$-matrix of $U_q(\su_2)$. The generators of the centralizer of $U_q(\su_2)$ in its three-fold tensor product 
are naturally expressed through conjugations of Casimir elements with $R$. They are seen as the images of the generators 
of $AW(3)$ under the embedding map by showing that they obey the $AW(3)$ relations. This is achieved by introducing a natural coaction also constructed with the help of the $R$-matrix.

\end{abstract}

\maketitle

%{\small MSC:\ 81R50;\ 81R10;\ 81U15.}
%

%{{\small  {\it \bf Keywords}: }}
%

\vspace{3mm}

\section{Introduction}
This letter addresses a long-standing question regarding the intrinsic description of the generators of a centrally extended Askey--Wilson algebra 
in its embedding into $U_q(\su_2)^{\otimes 3}$. The answer will be shown to involve Casimir elements and the universal $R$-matrix of $U_q(\su_2)$.

The Askey--Wilson algebra can be defined with three generators and relations. It has first been introduced \cite{Z91} as the algebra realized by the recurrence and 
$q$-difference operators intervening in the bispectral problem associated to the Askey--Wilson polynomials \cite{KLS}. 
This explains the name. Since the structure relations are not affected by truncations, this algebra also encodes the properties of the $q$-Racah polynomials. 
Owing to the connection with these 6$j$ or Racah coefficients for $U_q(\su_2)$ \cite{GZ}, a centrally extended Askey--Wilson algebra $AW(3)$  can be realized 
as the centralizer of the diagonal action of this quantum algebra in its three-fold tensor product. Related are the references \cite{koo1,koo2,ter,H-WH}. In this context, two generators of $AW(3)$ are naturally 
mapped under the coproduct onto the intermediate Casimir elements corresponding respectively to the recouplings of the first and last two factors in $U_q(\su_2)^{\otimes 3}$. 
A natural algebraic interpretation of the image of the third generator has however been lacking. This was circumvented so far by using  one of the relations which gives the third generator 
as the $q$-commutator of the other two; while this allows the homomorphism from $AW(3)$ into $U_q(\su_2)^{\otimes 3}$ to be defined, the resulting expression for this third generator is far from illuminating. 
Note that all three generators are needed to provide a $PBW$ basis for $AW(3)$. Besides the fact that this leaves a picture which is not fully satisfactory, this is a serious shortcoming 
in attempts to generalize $AW(3)$ to algebras of higher ranks. The natural approach - in fact the only one that has been conceived - is to define $AW(n)$ as the centralizer of $U_q(\su_2)$ in $U_q(\su_2)^{\otimes n}$. Proceeding with such an extension calls for an algebraic understanding of all centralizing elements in the tensor product. 
Significant progress towards describing the algebras $AW(n)$ have been made nevertheless. The algebra $AW(4)$ has been explored in \cite{PW} by including generators defined through the $q$-commutators of coproduct images of the Casimir element, as done for $AW(3)$, and obtaining from there various structure relations. 
Meaningful results have thus been found. The identification of the general $AW(n)$ has been attacked and largely advanced in \cite{DV,HD}. Much has been achieved in this case by cleverly designing a coaction map that has been used to define the generators, starting from the Casimir element of $U_q (\su_2)$, so as to ensure that these generators obey a $q$-deformation of natural structure relations 
(i.e. those of the generalized Bannai--Ito algebra $BI(n)$ \cite{DGV})  and by proving that this is so in many (but not all) cases. Still, in spite of this progress, an a priori algebraic description of the generators remained much desired.

We shall here settle this issue for $AW(3)$ by providing a simple expression for the image of its third generator in $U_q (\su_2)^{\otimes 3}$. 
The formula will involve conjugation with the universal $R$-matrix of $U_q(\su_2)$ and will be seen to explain the origin of the coaction introduced in \cite{DV}. 
Basic facts about $U_q(\su_2)$ and its universal $R$-matrix are collected in Section 1. Section 2 focuses on the centralizer of $U_q(\su_2)$  in $U_q(\su_2)^{\otimes 3}$; 
it provides the algebraic description that was missing. An additional centralizing element, conjugated to the usual third generator of $AW(3)$ is also identified; this will 
be related to observations made in \cite{PW}. The universal $R$-matrix and the Yang--Baxter equation are central here. With the expressions for the generators (in $U_q(\su_2)^{\otimes 3}$) in hand, 
Section 3 looks at their products and recovers the $AW(3)$ relations. To that end, a map from $U_q(\su_2)$  in $U_q(\su_2)^{\otimes 3}$ defined in terms of the $R$-matrix is introduced. 
It is pointed out that this map, once spelled out, coincides with the coaction used in \cite{DV}. The letter concludes with final remarks stressing the advantages of bringing the universal $R$-matrix 
in the description of the algebras $AW(n)$. As an illustration it is shown that a computation in $AW(4)$ can be performed with these tools in a comparatively much simpler way than otherwise.

\section{$U_q(\su_2)$ and its universal $R$-matrix}

In this section, we recall the definitions of the quantum algebra $U_q(\su_2)$ and of its universal $R$-matrix as well as some of their properties.
This allows to fix the notations and to make this letter more self-contained.

The associative algebra $U_q(\su_2)$ is generated by $E$, $F$ and $q^{H}$ with the defining relations
\begin{equation}
 q^{H}E=q  Eq^H\quad , \qquad q^{H}F=q^{-1} Fq^H  \quad \text{and} \qquad [E,F]=[2H]_q \ ,
\end{equation}
where $[X]_q=\frac{q^X-q^{-X}}{q-q^{-1}}$ and $q\neq \pm 1,\pm i$. The center of this algebra is generated by the following Casimir element
\begin{equation}
 C=-\frac{(q-q^{-1})^2}{q+q^{-1}}\left(FE + \frac{qq^{2H}+q^{-1}q^{-2H}}{(q-q^{-1})^2}  \right) \ .
 \label{eq:C} 
\end{equation}
The normalization of the Casimir element $C$ is irrelevant but chosen to yield computational simplifications.
There exists a homomorphism $\Delta:U_q(\su_2) \rightarrow U_q(\su_2)\otimes U_q(\su_2) $, called comultiplication, defined by 
\begin{equation}
 \Delta(E)=E \otimes q^{-H} +q^H \otimes E \quad , \qquad  \Delta(F)=F \otimes q^{-H} +q^{H} \otimes F \quad \text{and} \qquad \Delta(q^H)=q^H\otimes q^H\ .
\end{equation}
We recall that this comultiplication is coassociative 
\begin{equation}
 (\Delta \otimes \id)\Delta= (\id \otimes \Delta)\Delta \ .\label{eq:com}
\end{equation}

The quantum algebra $U_q(\su_2)$ is quasi-triangular because there exists a universal $R$-matrix $\cR\in U_q(\su_2) \otimes U_q(\su_2) $ which is invertible and satisfies 
\begin{eqnarray}
 \Delta(x) \cR  = \cR \Delta^{op}(x), \qquad \text{for } x\in U_q(\su_2)\ ,\label{eq:RD}
\end{eqnarray}
where the opposite comultiplication $\Delta^{op}(x)=x_{(2)} \otimes x_{(1)}$ if  $\Delta(x)=x_{(1)} \otimes x_{(2)}$ in the Sweedler notation, 
and
\begin{eqnarray}
 (\id \otimes \Delta)\cR=\cR_{12}\cR_{13}\ \quad \text{and}\qquad (\Delta \otimes \id)\cR=\cR_{23}\cR_{13}\ .\label{eq:idDRa}
\end{eqnarray}
In the previous relation \eqref{eq:idDRa}, we have used the usual notations $\cR_{12}=\cR^\alpha\otimes \cR_\alpha \otimes 1$, $\cR_{23}=1\otimes \cR^\alpha\otimes \cR_\alpha$
and $\cR_{13}=\cR^\alpha\otimes1\otimes  \cR_\alpha$ where $\cR=\cR^\alpha\otimes \cR_\alpha$ (the sum w.r.t. $\alpha$ is understood).
We will also use the following element
\begin{equation}
 \widetilde \cR= \cR_{21}=\cR_\alpha\otimes \cR^\alpha\ ,
\end{equation}
satisfying
\begin{eqnarray}
 \Delta^{op}(x) \widetilde \cR  = \widetilde \cR \Delta(x), \qquad \text{for } x\in U_q(\su_2)\ .\label{eq:RDop}
\end{eqnarray}
%One deduces that 
%\begin{eqnarray}
% (\id \otimes \Delta)\widetilde \cR=\widetilde \cR_{13}\widetilde \cR_{12}\ \quad \text{and}\qquad 
% (\Delta \otimes \id)\widetilde\cR=\widetilde\cR_{13}\widetilde\cR_{23}\ . \label{eq:idDR}
%\end{eqnarray}
The universal $R$-matrix also satisfies the Yang--Baxter equation
\begin{eqnarray}
 \cR_{12}\cR_{13}\cR_{23}&=&\cR_{23}\cR_{13}\cR_{12}\ , \label{YBE1}%\\
%  \cR_{13}\cR_{12}\widetilde \cR_{23}&=&\widetilde \cR_{23}\cR_{12}\cR_{13}\ , \label{YBE2}\\
% \widetilde\cR_{12}\cR_{23}\cR_{13}&=&\cR_{13}\cR_{23}\widetilde\cR_{12}\ . \label{YBE3}
\end{eqnarray}
% The universal R-matrix $\cR$ of $U_q(\su_2)$ 
and takes the following explicit form \cite{Dr}
\begin{equation}
 \cR=q^{2(H\otimes H)} \sum_{n=0}^\infty \frac{(q-q^{-1})^n}{[n]_q!} q^{n(n-1)/2} (E q^H  \otimes q^{-H} F)^n\ , \label{eq:uR}
\end{equation}
where $[n]_q!=[n]_q[n-1]_q\dots [2]_q[1]_q$ and, by convention, $[0]_q!=1$.
For future convenience, by using the commutation relations of $U_q(\su_2)$, we rewrite $\widetilde \cR$ as follows
\begin{equation}
 \widetilde \cR= \sum_{n=0}^\infty \frac{(q-q^{-1})^n}{[n]_q!} q^{n(n-1)/2} ( F q^{H} \otimes  q^{-H}  E  )^n\ \ q ^{2(H\otimes H)}= \Theta \ q ^{2(H\otimes H)}\ .\label{eq:Rt2}
\end{equation}

\section{Centralizer of $U_q(\su_2)$  in $U_q(\su_2)^{\otimes 3}$}

In this section, we want to describe the centralizer $\mC_3$ of the diagonal action of $U_q(\su_2)$ in $U_q(\su_2)^{\otimes 3}$:
\begin{equation}
 \mC_3=\{ X \in U_q(\su_2)^{\otimes 3}\ \big| \ [(\Delta\otimes \id)\Delta(x), X]=0\ , \ \  \forall x\in U_q(\su_2) \} . \label{eq:defc}
 \end{equation}
Let us define the so-called intermediate Casimir elements (in Sweedler's notation)
\begin{eqnarray}
 C_{12}=\Delta(C) \otimes 1=C_{(1)} \otimes C_{(2)} \otimes 1   \quad \text{and}\qquad C_{23}=1 \otimes \Delta(C) =1 \otimes C_{(1)} \otimes C_{(2)} \ ,
\end{eqnarray}
and the total Casimir element
\begin{eqnarray}
 C_{123}=(\Delta\otimes\id)\Delta(C)\ .\label{eq:Cr2}
\end{eqnarray}
We define also $C_1=C \otimes 1 \otimes 1$, $C_2=1 \otimes C \otimes 1$ and $C_3=1 \otimes 1 \otimes C$.
By using that the Casimir element $C$ is central in $U_q(\su_2)$, we deduce for example that
\begin{equation}
 [( \Delta\otimes \id)\Delta(x),C_{12}]=0 \quad\text{and}\qquad [( \id \otimes \Delta)\Delta(x),C_{23}]=0 \ , \ \ \forall x\in U_q(\su_2) \ .\label{eq:c12}
\end{equation}
 By definition \eqref{eq:defc}, $C_1$, $C_2$, $C_3$, $C_{12}$, $C_{23}$ and $C_{123}$ belong to the centralizer $\mC_3$
 with $C_1$, $C_2$, $C_3$ and $C_{123}$  belonging to the center of $\mC_3$. 
 It is well-known that these elements satisfy the Askey--Wilson algebra \cite{Z91}. We will come back to this point 
 in Section \ref{sec:AW}.\\

In the case of $U(\su_2)$ (\textit{i.e.} the limit $q\rightarrow 1$ of the case studied here),
 one can also prove that the intermediate Casimir $C_{13}= C_{(1)}\otimes 1 \otimes C_{(2)}$ belongs to the centralizer.
 For $q\neq 1$, it is not the case and the main objective of this letter is to provide a definition of this element for the quantum algebra.
 \begin{thm} The following elements of $U_q(\su_2)^{\otimes 3}$
 \begin{eqnarray}
  C_{13}^{(0)}&=&\widetilde \cR_{23}^{-1} C_{13} \widetilde \cR_{23}=\cR_{12}  C_{13} \cR_{12}^{-1}\ ,\label{eq:R131}\\
  C_{13}^{(1)}&=&\widetilde \cR_{12}^{-1} C_{13}\widetilde \cR_{12} =\cR_{23}  C_{13} \cR_{23}^{-1}\ .\label{eq:R132}
 \end{eqnarray}
 are in the centralizer $\mC_3$, where $\cR$ and $\widetilde \cR$ are defined in \eqref{eq:uR} and \eqref{eq:Rt2} and $C_{13}= C_{(1)}\otimes 1 \otimes C_{(2)}$.
 \end{thm}
\proof By using the coassociativity of the comultiplication \eqref{eq:com} and by conjugating with $\cR_{23}$, the first relation in \eqref{eq:c12} reads
\begin{equation}
 [( \id \otimes \Delta^{op})\Delta(x),\cR_{23}^{-1} C_{12 }\cR_{23}]=0\ .
\end{equation}
Finally, by exchanging the spaces $2$ and $3$, one gets that $C_{13}^{(0)}$ is in the centralizer
\begin{equation}
 [( \id \otimes \Delta)\Delta(x),\underbrace{\widetilde \cR_{23}^{-1} C_{13}\widetilde \cR_{23}}_{=C^{(0)}_{13}}]=0\ .
\end{equation}  
One proves similarly that $\cR_{12}  C_{13} \cR_{12}^{-1}$, $\widetilde \cR_{12}^{-1} C_{13}\widetilde \cR_{12}$ and $\cR_{23}  C_{13} \cR_{23}^{-1}$
are in the centralizer $\mC_3$. 

We must prove also the equality between $\widetilde \cR_{23}^{-1} C_{13} \widetilde \cR_{23}$ and $\cR_{12}  C_{13} \cR_{12}^{-1}$.
One gets
\begin{eqnarray}
 C_{13}^{(0)}=\widetilde \cR_{23}^{-1} C_{13}\widetilde \cR_{23}=\widetilde \cR_{23}^{-1} \left(C_{(1)}\otimes 1 \otimes C_{(2)}\right)\widetilde \cR_{23}=
 \widetilde \cR_{23}^{-1}\cR_{13} \left(C_{(2)}\otimes 1 \otimes C_{(1)}\right) \cR_{13}^{-1} \widetilde \cR_{23}\ ,
\end{eqnarray}
where we have used property \eqref{eq:RD}.
The Yang--Baxter equation \eqref{YBE1} implies that 
\begin{equation}
 C_{13}^{(0)}= \cR_{12}\cR_{13}\widetilde \cR_{23}^{-1} \cR_{12}^{-1} \left(C_{(2)}\otimes 1 \otimes C_{(1)}\right) \cR_{12}\widetilde \cR_{23}  \cR_{13}^{-1} \cR_{12}^{-1} \ .
\end{equation}
Now, from \eqref{eq:idDRa}, one deduces that $[\Delta(C)\otimes 1, (\Delta\otimes \id)(\cR)]=[\Delta(C)\otimes 1, \cR_{23}\cR_{13}]=0$ and that 
$[\left(C_{(2)}\otimes 1 \otimes C_{(1)}\right), \cR_{12}\widetilde \cR_{23}]=0$. Then, one obtains
\begin{equation}
 C_{13}^{(0)}= \cR_{12}\cR_{13} \left(C_{(2)}\otimes 1 \otimes C_{(1)}\right) \cR_{13}^{-1} \cR_{12}^{-1}=\cR_{12}  C_{13} \cR_{12}^{-1} \ .
\end{equation}
The equality between $\widetilde \cR_{12}^{-1} C_{13}\widetilde \cR_{12}$ and $\cR_{23}  C_{13} \cR_{23}^{-1}$ is proven similarly.
\endproof

%Using the previous theorem, it is easy to show that
%\begin{equation}
%	[C_{123},C_{13}^{(0)}]=0 \qquad \text{and} \quad [C_{123},C_{13}^{(1)}]=0 \ .
%\end{equation}
From relations \eqref{eq:R131} and \eqref{eq:R132}, one deduces that $C_{13}^{(0)}$ and $C_{13}^{(1)}$ are conjugated:
\begin{equation}
 C_{13}^{(1)}= \cR_{23}\widetilde \cR_{23} C_{13}^{(0)} (\cR_{23}\widetilde \cR_{23})^{-1}=(\cR_{12}\widetilde \cR_{12})^{-1} C_{13}^{(0)}  \cR_{12}\widetilde \cR_{12} \ .
\end{equation}

%It is very natural to introduce the element $\cX$ since $\cX_{12}$ and $\cX_{23}$ are in the centralizer $\mC_3$. It is easy demonstrated by using \eqref{eq:RD} and \eqref{eq:RDop}.

%\begin{prop}
% The following elements, for $n=1,2,\dots$
% \begin{eqnarray}
%  C^{(n)}_{13}&=&\cX_{23}^{n-1} C_{13}^{(1)}  \cX_{23}^{-n+1}\\
%  C^{(n)}_{23}&=&\cX_{12}^{-n+1} C_{23}^{(1)} \cX_{12}^{n-1}\\
%  C^{(n)}_{12}&=&\cX_{23}^{n-1} C_{12}^{(1)}  \cX_{23}^{-n+1}
% \end{eqnarray}
% where, by convention $C_{12}^{(1)}=C_{12}$ and $C_{23}^{(1)}=C_{23}$, are in the centralizer $\mC_3$.
%\end{prop}
%\proof  ?? A reflechir si ce sont les seuls elements ???
%\endproof

\section{The Askey--Wilson algebra $AW(3)$ \label{sec:AW} }

In this section, we study the algebra satisfied by the intermediate Casimir elements introduced in the previous section and 
connect it with the central extension $AW(3)$ of the Askey--Wilson algebra introduced in  \cite{Z91}. We start by proving the following lemma.
\begin{lem}\label{lem1}
 The map defined by
 \begin{eqnarray}\label{eq:hatt}
  \tau\ :\ U_q(\su_2)& \rightarrow & U_q(\su_2)\otimes U_q(\su_2)\\
  x&\mapsto &\widetilde \cR^{-1} (1 \otimes x) \widetilde \cR \nonumber
 \end{eqnarray}
 yields the following explicit expressions when acting on the different elements of $U_q(\su_2)$ listed below:
 \begin{eqnarray}
 \tau(C)&=& 1\otimes C\ ,\label{eq:tau1}\\
  \tau(q^{-H}E)&=&q^{-2H} \otimes q^{-H} E\ ,\label{eq:tau2}\\
  \tau(q^{-2H})&=& 1 \otimes q^{-2H} -(q-q^{-1})^2  q^{-H}F \otimes  q^{-H}E \ ,\label{eq:tau3}\\
  \tau(Fq^{-H})&=& q^{2H}\otimes Fq^{-H} + q^{-1}(q+q^{-1})Fq^H\otimes(C+q^{-2H})-(q-q^{-1})^2F^2\otimes q^{-H}E \ .\label{eq:tau4}
\end{eqnarray}
\end{lem}
\proof We must prove that the map given in the theorem reproduces relations \eqref{eq:tau1}-\eqref{eq:tau4}. For relation \eqref{eq:tau1},
it is direct, knowing that $C$ commutes with any element of $U_q(\su_2)$. To prove relation \eqref{eq:tau2}, one computes (using the explicit
form \eqref{eq:Rt2} of $\widetilde \cR$)
\begin{eqnarray}
  \tau(q^{-H}E)= \widetilde \cR^{-1} (1 \otimes q^{-H}E )\Theta \ q ^{2(H\otimes H)}
 = \widetilde \cR^{-1}\Theta \  (1 \otimes q^{-H}E) q ^{2(H\otimes H)}=q^{-2H}\otimes q^{-H}E\ ,
\end{eqnarray}
which reproduces \eqref{eq:tau2}.

We want now to compute $\tau(q^{-2H})$:
\begin{eqnarray}
\tau(q^{-2H})&=&\widetilde \cR^{-1} (1 \otimes q^{-2H} )\widetilde\cR= \widetilde \cR^{-1} (1 \otimes q^{-2H} )q^{2(H\otimes H)} \sum_{n=0}^\infty a_n (q^{-H} F \otimes  E q^H  )^n\\
&=&\widetilde \cR^{-1} q^{2(H\otimes H)} \sum_{n=0}^\infty a_n q^{-2n} (q^{-H} F \otimes  E q^H  )^n\  (1 \otimes q^{-2H} )\ ,
\end{eqnarray}
where we have introduced the parameters
\begin{equation}
a_n = \frac{(q-q^{-1})^n}{[n]_q!} q^{n(n-1)/2}\ .
\end{equation}
Remarking that 
\begin{equation}
	a_n q^{-2n}=a_n-a_n[n]_q q^{-n}(q-q^{-1}) \ ,
	\label{eq:rel1_an}
\end{equation}
one gets
\begin{equation}
\tau(q^{-2H})=\widetilde \cR^{-1}\left( \widetilde \cR- q^{2(H\otimes H)} \sum_{n=0}^\infty a_{n+1}[n+1]_q q^{-(n+1)}(q-q^{-1}) (q^{-H} F \otimes  E q^H  )^{n+1} \right)  (1 \otimes q^{-2H} )\ .
\end{equation}
It is easy to show that the parameters $a_n$ satisfy $a_{n+1}[n+1]_q = q^n(q-q^{-1})a_n$, which allows to recover \eqref{eq:tau3}.

Similarly, to prove \eqref{eq:tau4}, one computes
\begin{eqnarray}
\tau(Fq^{-H})&=&\widetilde \cR^{-1}(1\otimes Fq^{-H})\widetilde \cR=\widetilde \cR^{-1} (1\otimes Fq^{-H}) q^{2(H\otimes H)} \sum_{n=0}^\infty a_n (q^{-H} F \otimes  E q^H  )^n\\
&=&\widetilde \cR^{-1} q^{2(H\otimes H)} (q^{2H}\otimes Fq^{-H}) \sum_{n=0}^\infty a_n (q^{-H} F \otimes  E q^H  )^n \ .
\end{eqnarray}
Then, the identity
\begin{eqnarray}
[F,E^n]=\frac{[n]_q}{q-q^{-1}}(q^{n-1}q^{-2H}-q^{-(n-1)}q^{2H})E^{n-1}
\end{eqnarray}
can be used to write
\begin{equation}
\tau(Fq^{-H})=\widetilde \cR^{-1} q^{2(H\otimes H)} \sum_{n=0}^\infty a_{n}(q^{-H} F \otimes  E q^H  )^n\left( q^{-2n}q^{2H}\otimes Fq^{-H}+q^{-2} Fq^H\otimes (q^{-2n}q^{-2H}-q^{2H})\right) \ .
\end{equation}
Using again relation \eqref{eq:rel1_an}, one finds 
\begin{equation}
\tau(Fq^{-H})=q^{2H}\otimes Fq^{-H}-q^{-1}(q-q^{-1})^2 Fq^{H} \otimes FE-Fq^H\otimes (q^{2H}-q^{-2H})-(q-q^{-1})^2F^2\otimes q^{-H}E \ .
\end{equation}
Finally, expressing $FE$ in terms of $C$ from definition \eqref{eq:C}, one recovers \eqref{eq:tau4}.  
\endproof
Using Lemma \ref{lem1}, we can rewrite $C_{13}^{(0)}$ \eqref{eq:R131} as follows
\begin{eqnarray}
 &&C_{13}^{(0)}=(1 \otimes \tau ) \Delta(C)\ \label{eq:C13}\\
 &=& (q^{2H}+C)\otimes \tau(q^{-2H})+q^{2H} \otimes \tau(C)-\frac{(q-q^{-1})^2}{q+q^{-1}}\left(q^HE\otimes\tau(Fq^{-H})+Fq^H\otimes\tau(q^{-H}E)\right).\label{eq:C13exp}
\end{eqnarray}

\begin{prop} The following relation
 \begin{eqnarray}
 &&\frac{1}{q-q^{-1}}  [C_{12},C_{23}]_q = C_{13}^{(0)} +C_1C_3+C_2C_{123} \label{eq:AW31}
\end{eqnarray}
holds in  $U_q(\su_2)^{\otimes 3}$.
\end{prop}
\proof Using the expressions for the maps under $\tau$ given in Lemma \ref{lem1}, we obtain $C_{13}^{(0)}$ in terms of the generators of $U_q(\su_2)$.
A direct computation using the commutation relations of $U_q(\su_2)$ proves the relation of the proposition.\endproof

One of the advantages of the construction with the universal $R$-matrix is that we can deduce all the other defining relations of $AW(3)$ from \eqref{eq:AW31} and some other relations.
\begin{corollary}  The following relations
 \begin{eqnarray}
 && \frac{1}{q-q^{-1}} [C_{13}^{(0)},C_{12}]_q= C_{23} +C_2C_3+C_1C_{123},\label{eq:AW32}\\
  &&\frac{1}{q-q^{-1}} [C_{23},C_{13}^{(0)}]_q = C_{12} +C_1C_2+C_3C_{123},\label{eq:AW33}\\
 &&\frac{1}{q-q^{-1}}  [C_{23},C_{12}]_q = C_{13}^{(1)} +C_1C_3+C_2C_{123}, \label{eq:AW41}\\
  &&\frac{1}{q-q^{-1}} [C_{12},C_{13}^{(1)}]_q= C_{23} +C_2C_3+C_1C_{123},\label{eq:AWI32}\\  
  &&\frac{1}{q-q^{-1}} [C_{13}^{(1)},C_{23}]_q= C_{12} +C_1C_2+C_3C_{123},\label{eq:AWI33}
\end{eqnarray}
hold in  $U_q(\su_2)^{\otimes 3}$.
\end{corollary}
\proof 
We use the second relation in \eqref{eq:R131} as well as the definitions \eqref{eq:Cr2} to write relation \eqref{eq:AW31} as follows
\begin{equation}
 \frac{1}{q-q^{-1}}  [\Delta(C)\otimes 1,C_{23}]_q = \cR_{12}  C_{13} \cR_{12}^{-1}+C_1C_3+C_2(\Delta \otimes \id)\Delta(C)\ .
\end{equation}
Exchanging the spaces 1 and 2, the previous relation becomes
\begin{equation}
\frac{1}{q-q^{-1}}  [\Delta^{op}(C)\otimes 1,C_{13}]_q = \widetilde \cR_{12}  C_{23} \widetilde \cR_{12}^{-1}+C_2C_3+C_1(\Delta^{op} \otimes \id)\Delta(C) \ ,
\end{equation}
which leads to \eqref{eq:AWI32} after conjugating by $\widetilde\cR_{12}$ (using property \eqref{eq:RDop}).

Then, one starts from the relation \eqref{eq:AWI32} we have just proven, uses the second relation in \eqref{eq:R132} to express $C_{13}^{(1)}$ and exchanges spaces 2 and 3 to write
\begin{equation}
	\frac{1}{q-q^{-1}} [ C_{13}, \widetilde R_{23}C_{12}\widetilde R_{23}^{-1}]_q= 1\otimes \Delta^{op}(C) +C_2C_3+C_1(\id \otimes \Delta^{op})\Delta(C) \ .
\end{equation}
Conjugating with $\widetilde R_{23}$, one proves relation \eqref{eq:AW32}.
Performing the same two steps starting from \eqref{eq:AW32}, one proves \eqref{eq:AW41} and \eqref{eq:AW33}. Finally, the two same steps prove \eqref{eq:AWI33} and give again the equation \eqref{eq:AW31}.
\endproof

We now have a number of remarks regarding the merits of the $R$-matrix approach developed above.
\begin{rem} Relations \eqref{eq:AW31}, \eqref{eq:AW32} and \eqref{eq:AW33} are the defining relations of central extension $AW(3)$ 
of the Askey--Wilson algebra introduced in \cite{Z91}.
Therefore, the results presented in this letter offer another proof that the intermediate Casimir elements of $U_q(\su_2)$ provide a realization of $AW(3)$. 
In previous works \cite{Z91,ter, DV, DD, PW}, $C_{13}^{(0)}$ was defined by relation \eqref{eq:AW31} whereas in our approach, it is defined independently of the commutation relations
via relation \eqref{eq:R131}.
\end{rem}

\begin{rem}\label{rem2} The map $\tau$ with images given by \eqref{eq:tau1}-\eqref{eq:tau4} has in fact been introduced in \cite{DV,DD} so as to obtain $C_{13}^{(0)}$ as in relation \eqref{eq:C13}.  
Our definition \eqref{eq:hatt} gives a nice and powerful interpretation of this map. Let us remark that the comultiplication used in this letter is slightly different from the ones used in \cite{DV,DD}. 
In order to establish exactly the correspondence, the following transformation on our generators and Casimir element must be performed:
$q^{2H}\rightarrow K$, $E\rightarrow EK^{-1/2}$, $F\rightarrow K^{1/2}F$, $q\rightarrow Q$ and
$C\rightarrow -\Lambda/(Q+Q^{-1})$.
\end{rem}

\begin{rem} To illustrate the appropriateness and advantages of definition \eqref{eq:hatt} of the map $\tau$, we here prove its coaction property in a much simpler way than the direct calculation described in \cite{DV,DD}.
 Using relation \eqref{eq:idDRa}, it is easy to compute, for $x\in U_q(\su_2)$,
\begin{equation}
(\Delta \otimes \id) \tau(x)=(\Delta \otimes \id)\left(\widetilde\cR^{-1} (1 \otimes x) \widetilde \cR\right)
=\widetilde\cR_{23}^{-1}\widetilde\cR_{13}^{-1}(1 \otimes 1 \otimes x)\widetilde\cR_{13}\widetilde\cR_{23}\ ,
\end{equation}
and 
\begin{equation}
( \id \otimes  \tau ) \tau(x)=( \id \otimes \tau)\left(\widetilde\cR^{-1} (1 \otimes x) \widetilde \cR\right)
=\widetilde\cR_{23}^{-1}\widetilde\cR_{13}^{-1}(1 \otimes 1 \otimes x)\widetilde\cR_{13}\widetilde\cR_{23}\ .
\end{equation}
This proves that $(\Delta \otimes \id) \tau(x)=( \id \otimes  \tau ) \tau(x)$ and thus that $\tau$ is a left coaction.
\end{rem}

\begin{rem}
We can define also a right coaction $\check \tau$ given by
\begin{eqnarray}
\check \tau\ :\ U_q(\su(2))& \rightarrow & U_q(\su(2))\otimes U_q(\su(2))\\
x&\mapsto & \cR (x\otimes 1) \cR^{-1} \nonumber \ ,
\end{eqnarray}
satisfying
\begin{equation}
(\check  \tau \otimes \id  ) \check  \tau = (\id \otimes \Delta ) \check  \tau \ . 
\label{eq:coactiond}
\end{equation}
We can show following steps similar to those of the proof of Lemma \ref{lem1} that this right coaction coincides with the one introduced recently in \cite{HD}
with the identifications: $q^{2H}\rightarrow K$, $E\rightarrow EK^{-1/2}$, $F\rightarrow K^{1/2}F$ and
$C\rightarrow -\Lambda/(q+q^{-1})$.
\end{rem}

\begin{rem}
 The element $C_{13}^{(1)}$ has been introduced previously in \cite{PW} (where it is called $IQ^{(13)}$) and defined by relation \eqref{eq:AW41}.
 Our definition \eqref{eq:R132} gives a new interpretation of this element.
\end{rem}

\section{Conclusion and perspective}

In this letter, we study the centralizer of the diagonal action of $U_q(\su_2)$ and its connection with the Askey--Wilson algebra $AW(3)$.
In comparison with the previous approaches, we have emphasized the relevance of the universal $R$-matrix of $U_q(\su_2)$. We believe that its use
offers a deeper understanding of the realization of the Askey--Wilson algebra in terms of the intermediate Casimir elements.
It should moreover simplify the computations for further investigations. 
To illustrate this point, let us show how one computation can be simplified with this approach in the higher rank generalization $AW(4)$ of $AW(3)$ 
examined in \cite{PW}. 
The algebra $AW(4)$ can be embedded in $U_q(\su_2)^{\otimes 4}$ and, in particular, one defines
\begin{eqnarray}
  C_{13}^{(0)}&=&\widetilde \cR_{23}^{-1} C_{13} \widetilde \cR_{23}=\cR_{12}  C_{13} \cR_{12}^{-1}\ ,\label{eq:R131_4}\\
  C_{24}^{(1)}&=&\widetilde \cR_{23}^{-1} C_{24}\widetilde \cR_{23} =\cR_{34}  C_{24} \cR_{34}^{-1}\ .\label{eq:R132_4}
 \end{eqnarray}
Looking at the commutation relations, we can prove that these elements correspond to $Q^{(13)}$ and $IQ^{(24)}$ of \cite{PW}.
In the formalism introduced here, we see immediately that 
\begin{equation}
  [C_{13}^{(0)},C_{24}^{(1)}]=0\ ,
 \end{equation}
whereas the proof without the use of the $R$-matrix presented in \cite{PW} is quite cumbersome. We believe that the $R$-matrix approach 
we have elaborated will prove quite helpful in the study of the higher rank generalizations of $AW(3)$. In a related series of papers \cite{CPV}, \cite{CFV}, 
the Temperley--Lieb algebra with $q=1$, the Brauer algebra (and others) over $3$ strands have been identified as  quotients of the Racah \cite{GVZ2} 
and Bannai--Ito \cite{TVZ} algebras of rank 1. The results reported here pave the way to the pursuit of this program for $AW(3)$ as well as in situations of higher ranks with an arbitrary number of strands. It is our intent to actively continue this research.
Let us mention finally that, in a companion letter \cite{CVZ}, we have provided a parallel description of the Bannai--Ito algebras using the universal $R$-matrix of the Lie superalgebra $\mathfrak{osp}(1|2)$.

\bigskip 

{\bf Acknowledgments:}
We have much benefited from discussions with L. Frappat and E. Ragoucy.
N.~Cramp\'e is gratefully holding a CRM--Simons professorship.
The research of L.~Vinet is supported
in part by a Discovery Grant from the Natural Science and Engineering
Research Council (NSERC) of Canada. J.~Gaboriaud and M.~Zaimi hold a NSERC graduate scholarship.

\end{document}